\documentclass{amsart}
\usepackage{amsmath}
\usepackage{amsfonts}
\usepackage{amssymb}
\usepackage{color}
\usepackage{graphics,graphicx}
\usepackage{hyperref}
\usepackage{amssymb,latexsym,amsthm,xypic}

  \usepackage{mathrsfs}

\newtheorem{thm}{Theorem}[section]
\newtheorem{cor}[thm]{Corollary}
\newtheorem{lem}[thm]{Lemma}
\newtheorem{prop}[thm]{Proposition}

\theoremstyle{definition}
\newtheorem{defn}{Definition}[section]

\theoremstyle{remark}
\newtheorem{rem}{Remark}[section]

\newcommand{\be}{\begin{equation}}
\newcommand{\ee}{\end{equation}}
\newcommand{\bea}{\begin{eqnarray}}
\newcommand{\eea}{\end{eqnarray}}
\newcommand{\ben}{\begin{eqnarray*}}
	\newcommand{\een}{\end{eqnarray*}}
\newcommand{\bt}{\begin{split}}
	\newcommand{\et}{\end{split}}
\newcommand{\bet}{\begin{equation}}
\newcommand{\mc}{\mathbb{C}}
\newcommand{\mr}{\mathbb{R}}
\newcommand{\ra}{\rightarrow}

\begin{document}

\title[Linear invariants and plurisubharmonic variation]{Linear invariants of complex manifolds and their  plurisubharmonic  variations}
\author[F. Deng]{Fusheng Deng}
\address{Fusheng Deng: \ School of Mathematical Sciences, University of Chinese Academy of Sciences\\ Beijing 100049, P. R. China}
\email{fshdeng@ucas.ac.cn}
\author[Z. Wang]{Zhiwei Wang}
\address{ Zhiwei Wang: \ School
	of Mathematical Sciences\\Beijing Normal University\\Beijing\\ 100875\\ P. R. China}
\email{zhiwei@bnu.edu.cn}
\author[L. Zhang]{Liyou Zhang}
\address{ Liyou Zhang: \ School of Mathematical Sciences\\Capital Normal University\\Beijing \\100048\\
	P. R.  China}
\email{zhangly@cnu.edu.cn}
\author[X. Zhou]{Xiangyu Zhou}
\address{Xiangyu Zhou:\ Institute of Mathematics, AMSS, and Hua Loo-Keng Key Laboratory of Mathematics, Chinese Academy of Sciences, Beijing 100190, China}
\email{xyzhou@math.ac.cn}

\begin{abstract}
For a bounded domain $D$ and a real number $p>0$,
we denote by $A^p(D)$ the space of $L^p$ integrable holomorphic functions on $D$, equipped with the $L^p$- pseudonorm.
We prove that two bounded hyperconvex domains $D_1\subset \mc^n$ and $D_2\subset \mc^m$ are biholomorphic
(in particular $n=m$) if there is a linear isometry between $A^p(D_1)$ and $A^p(D_2)$ for some $0<p<2$.
The same result holds for $p>2, p\neq 2,4,\cdots$, provided that the $p$-Bergman kernels on $D_1$ and $D_2$ are exhaustive.
With this as a motivation, we show that, for all $p>0$, the $p$-Bergman kernel on a strongly pseudoconvex domain with $\mathcal C^2$ boundary
or a simply connected homogeneous regular domain is exhaustive.
These results shows that spaces of pluricanonical sections of complex manifolds equipped with canonical pseudonorms are important
invariants of complex manifolds.
The second part of the present work devotes to studying variations of these invariants.
We show that the direct image sheaf
of the twisted relative $m$-pluricanonical bundle associated to a holomorphic family of Stein manifolds or compact K\"ahler manifolds
is positively curved, with respect to the canonical singular Finsler metric.
\end{abstract}

\thanks{The author was partially supported by the Fundamental Research Funds for the Central Universities and by the NSFC grant NSFC-11701031}


\maketitle
\tableofcontents

\section{Introduction}
It is well known that two pseudoconvex domains (or even Stein manifolds) $\Omega_1$ and $\Omega_2$ are biholomorphic
if and only if $\mathcal O(\Omega_1)$ and $\mathcal O(\Omega_2)$, the spaces of holomorphic functions,
are isomorphic as $\mc$-algebras with unit.
This implies that the holomorphic structure of a pseudoconvex domain
is uniquely determined by algebraic structure of the space of holomorphic functions on it.
In the first part of the present work,
we prove some results in a related but different direction.

To state the results, we need to recall some notions.
Let $\Omega$ be a bounded domain in $\mathbb C^n$.
Let $z=(z_1,\cdots,z_n)$ be the natural holomorphic coordinates of $\mathbb C^n$,
and let $d\lambda_n:={(\frac{i}{2})^n}dz_1\wedge d\overline{z}_1\wedge\cdots\wedge dz_n\wedge d\overline{z}_n$
be the canonical volume form on $\Omega$.
For $p>0$,  denote by  $A^p(\Omega)$  the space of all holomorphic functions $\phi$ on $\Omega$ with finite $L^p$-norm
\begin{align*}
\|\phi\|_p:=\left(\int_\Omega |\phi|^pd\lambda_n\right)^{1/p}.
\end{align*}
It is a standard fact that for $p\geq 1$,  $A^p(\Omega)$ are seperable Banach spaces,
and for $0<p<1$, $A^p(\Omega)$ are complete separable metric spaces with respect to the metric
$d(\varphi_1,\varphi_2):=\|\varphi_1-\varphi_2\|_p^p$.

The $p$-Bergman kernel is defined as follows:
\begin{align*}
B_{\Omega,p}(z):=\sup_{\phi\in A^p(\Omega)}\frac{|\phi(z)|^2}{\|\phi\|^2_p}.
\end{align*}
When $p=2$, $B_{\Omega,p}$ is the ordinary Bergman kernel.
By a standard argument of Montel theorem,
one can prove that  $B_{\Omega,p}$ is a continuous plurisubharmonic function on $\Omega$.
We say that $B_{\Omega,p}$ is exhaustive if for any real number $A$ the set $\{z\in\Omega| B_{\Omega,p}(z)\leq A\}$ is compact.

A bounded domain $\Omega$ is called hyperconvex if there is a
plurisubharmonic function $\rho:\Omega\ra [-\infty, 0)$ such that
for any $c<0$ the set $\{z\in \Omega|\rho(z)\leq c\}$ is compact.

\begin{thm}\label{thm-intro: main theorem 1}
Let $\Omega_1\subset\mc^n$ and $\Omega_2\subset\mc^m$ be bounded hyperconvex domains.
Suppose that there is a $p>0, p\neq 2,4,6 \cdots$ , such that
\begin{itemize}
\item[(1)] there is a  linear isometry $T:A^p(\Omega_1)\ra A^p(\Omega_2)$,  and
\item[(2)] the $p$-Bergman kernels of  $\Omega_1$ and $\Omega_2$ are exhaustive,
\end{itemize}
then $m=n$ and there exists a unique biholomorphic map $F:\Omega_1\ra \Omega_2$ such that
\begin{align*}
|T\phi\circ F||J_F|^{2/p}=|\phi|, \ \forall \phi\in A^p(\Omega_1),
\end{align*}
where $J_F$ is the holomorphic Jacobian of $F$.
If $n=1$, the assumption of hyperconvexity can be dropped.
\end{thm}

Throughout this paper, that $T:A^p(\Omega_1)\ra A^p(\Omega_2), \phi\mapsto T\phi$ is a linear isometry
means that $T$ is a linear surjective map and $\|T\phi\|_p=\|\phi\|_p$ for all $\phi\in A^p(\Omega_1)$.

In \S \ref{subset:define map and equidim},
we will show that for any two bounded domains $\Omega_1$ and $\Omega_2$,
$dim \Omega_1=dim \Omega_2$ provided that $A^p(\Omega_1)$ and $A^p(\Omega_2)$
are linear isometric for some $p>0, p\neq 2,4, \cdots$.

Theorem \ref{thm-intro: main theorem 1} implies that the exhaustion of the $p$-Bergman kernels
for bounded domains is a very important property.
For the case that $p<2$, Ning-Zhang-Zhou get a complete result by showing that
a bounded domain $\Omega$ is pseudoconvex if and only if its $p$-Bergman kernel
is exhaustive \cite{NZZ16}, the proof of which is based on a $L^p$-variant of
the Ohsawa-Takegoshi extension theorem.
Therefore we get from \ref{thm-intro: main theorem 1} the following

\begin{thm}\label{thm-intro:p<2 case}
Let $\Omega_1\subset\mc^n$ and $\Omega_2\subset\mc^m$ be bounded hyperconvex domains.
Suppose that there is a linear isometry $T:A^p(\Omega_1)\ra A^p(\Omega_2)$ for some $p\in(0,2)$ ,
then $m=n$ and there exists a unique biholomorphic map $F:\Omega_1\ra \Omega_2$ such that
\begin{align*}
|T\phi\circ F||J_F|^{2/p}=|\phi|, \ \forall \phi\in A^p(\Omega_1).
\end{align*}
If $n=1$, the assumption of hyperconvexity can be dropped.
\end{thm}

It is known that bounded pseudoconvex domains with Lipschitz boundary
are hyperconvex \cite{Demailly87}.
Even more, bounded pseudoconvex domains with $\alpha$-H\"order boundary for all $\alpha<1$ are hyperconvex \cite{Ave-Hed-Per15}.
Recall that a bounded domain $\Omega\subset \mc^n$ is called \emph{homogeneous regular},
a concept proposed by Liu-Sun-Yau \cite{Liu-Sun-Yau04},
if there is a constant $c\in (0,1)$ such that for any $z\in \Omega$,
there is a holomorphic injective map $f:\Omega\ra \mathbb B^n$ with
$f(z)=0$ and $\mathbb B^n(c)\subset f(\Omega)$,
where $\mathbb B^n(c)=\{z\in\mc^n; \|z\|<c\}$.
It is shown in \cite{Yeung09} that all homogeneous regular domains are hyperconvex.

The exhaustion of $p$-Bergman kernels for $p>2$ is a subtle problem.
In the present paper, we show that $p$-Bergman kernels of simply connected homogeneous regular domains
or bounded strongly pseudoconvex domains are exhaustive for all $p>0$.

\begin{thm}\label{thm-intro: exhaustion simply connected USD}
If $\Omega\subset\mc^n$ is a simply connected and homogeneous regular domain,
then the $p$-Bergman kernel of $\Omega$ is exhaustive.
 \end{thm}

It seems interesting to consider whether the condition about simply connectedness
in Theorem \ref{thm: exhaustion simply connected USD} can be dropped.
This is the case for bounded strongly pseudoconvex domains with $\mathcal C^2$ boundary, which are known to be homogeneous regular \cite{Deng-Guan-Zhang16}.
The proof is based on results on exposing boundary points of strongly pseudoconvex domains in \cite{Diedrich-Fornaess-Wold13}\cite{DFW18}.

\begin{thm}\label{thm-intro: exhaustion strongly pscd}
For any bounded strongly pseudoconvex domain $\Omega$ in $\mathbb C^n$
with $\mathcal C^2$ smooth boundary, the $p$-Bergman kernel is exhaustive for all $p>0$.
\end{thm}

We conjecture that Theorem \ref{thm-intro: exhaustion strongly pscd} holds for all bounded pseudoconvex domains 
with $\mathcal C^2$ smooth boundary.

We have known many homogeneous regular domains, such as
Teichm\"{u}ller spaces of compact Riemann surfaces,
bounded homogeneous domains,
strongly pseudoconvex domains with $\mathcal C^2$ boundary,
bounded convex domains,
and bounded complex convex domains(see \cite{Deng-Guan-Zhang16},\cite{Frankel91},\cite{Kim-Zhang16},\cite{Nikolov-Andreev16}).

By Theorem \ref{thm-intro: main theorem 1}, Theorem \ref{thm-intro: exhaustion simply connected USD},
and Theorem \ref{thm-intro: exhaustion strongly pscd}, we get

\begin{thm}\label{thm-intro:usd of spsc domain}
Let $\Omega_1\subset\mc^n$ and $\Omega_2\subset\mc^m$ be bounded domains
which are simply connected and homogeneous regular, or
strongly pseudoconvex with $\mathcal C^2$ smooth boundary.
Suppose that there is a linear isometry $T:A^p(\Omega_1)\ra A^p(\Omega_2)$ for some $p>0, p\neq 2,4, \cdots$ ,
then $m=n$ and there exists a unique biholomorphic map $F:\Omega_1\ra \Omega_2$ such that
\begin{align*}
|T\phi\circ F||J_F|^{2/p}=|\phi|, \ \forall \phi\in A^p(\Omega_1).
\end{align*}
\end{thm}

For the case that $n=m=1$ and $p=1$,
Theorem \ref{thm-intro:p<2 case} was proved by Lakic \cite{Lakic97} and Markovic in \cite{Mark03}.
Application of spaces of pluricanonical forms with pseudonorms to birational geometry of projective algebraic
manifolds  was proposed by Chi and Yau in \cite{Chi-Yau2008} and was studied in \cite{Chi2016}\cite{Yau15}.
These works form the motivation to Theorem \ref{thm-intro: main theorem 1}.

Our method to Theorem \ref{thm-intro: main theorem 1} is partially inspired by the argument in \cite{Mark03}.
The  argument in \cite{Mark03} can be divided into two main parts as follows.
The first part is to construct a biholomorphic map, say $F$, from
a dense open subset $\Omega'_1\subset \Omega_1$ onto a dense open subset $\Omega'_2\subset \Omega_2$;
and the second part is to show that $F$ extends to a biholomorphic map from the whole $\Omega_1$ to $\Omega_2$.
The first part is based on a result of Rudin on equimeasurability \cite{Rud76}.
The second part, namely continuing $F$ across $\Omega_1\backslash \Omega'_1$, is more technical
and involves the uniformization theorem for Riemann surfaces.

The proof of Theorem \ref{thm-intro: main theorem 1} in the present paper is also divided into
two parts in the same way.
The construction of $F$ in the first part is based on the method in \cite{Mark03}.
The difference is that our starting point is an intrinsic definition of $F$,
where that in \cite{Mark03} is a countable dense subset of $A^1(\Omega)$, from which the map $F$ is defined.
The benefit of the modification is that $F$ is obviously uniquely determined by $T$
and some important properties of $F$ become clear just from the definition.
The method to the second part, namely the part about extension of $F$,
is totally different from that in \cite{Mark03}.
The new observation comes from pluripotential theory.
Using the assumption of exhaustion of $p$-Bergman kernels,
we construct a plurisubharmonic functions say $\rho$ on $\Omega_1$ such that
$\rho=-\infty$ exactly on $\Omega_1\backslash \Omega'_1$,
and hence $\Omega_1\backslash \Omega'_1$ is a closed pluripolar set in $\Omega_1$.
Then, by a basic result in pluripotential theory about extension of holomorphic functions,
we can continue $F$ to a holomorphic map on the whole $\Omega_1$.
The hyperconvexity property is then used to prove that $\Omega'_j=\Omega_j,\ j=1,2$.

\begin{rem}
(1). Let $f:\Omega_1\ra\Omega_2$ be a biholomorphic map between two bounded domains.
   If $\Omega_1$ is simply connected or $p=2/m$ for some integer $m$.
   then $f$ induces a linear isometry $T:A^p(\Omega_2)\ra A^p(\Omega_1)$ as:
   $$\psi\mapsto (\psi\circ f)J^{2/p}_f,\ \psi\in A^p(\Omega_2).$$

(2). When $p=2/m$ for some integer $m$, $A^p(\Omega)$ can be defined intrinsically as the space of
   holomorphic sections of $mK_\Omega$ with finite intrinsic pseudonorm (will be defined later),
   where $K_\Omega$ is the canonical bundle of $\Omega$.

(3). For any two bounded domains $\Omega_1$ and $\Omega_2$,  $A^2(\Omega_1)$ and $A^2(\Omega_2)$
   are always linear isometric since they are both sparable Hilbert spaces.
\end{rem}

The second part of the present paper devotes to the study of variation of $A^p(\Omega)$ as $\Omega$ varies.
We will consider holomorphic families of more general complex manifolds.

We first consider a family of pseudoconvex domains.
Let $\Omega\subset \mathbb{C}^{r+n}=\mathbb{C}_t^r\times \mathbb{C}_z^n$ be a pseudoconvex domain.
Let $p:\Omega\rightarrow \mathbb{C}^r$ be the natural  projection.
We denote $p(\Omega)$ by $D$ and denote $p^{-1}(t)$ by $\Omega_t$ for $t\in D$.
Let $\varphi$ be a plurisubharmonic function on $\Omega$ and let $m\geq 1$ be a fixed integer.
For an open subset $U$ of $D$, we denote by $\mathcal F(U)$ the space of holomorphic functions $F$ on $p^{-1}(U)$ such that
$\int_{p^{-1}(K)}|F|^{2/m}e^{-\varphi}\leq\infty$ for all compact subset $K$ of $D$.
For $t\in D$, let
$$E_{m,t}=\{F|_{\Omega_t}:F\in\mathcal F(U),\ U\subset D\ \text{open\ and}\ U\ni t\}.$$
Then $E_{m,t}$ is a vector space with the following pseudonorm:
$$H(f):=|f|_m=\left(\int_{D_t}|f|^{2/m}e^{-\varphi_t}\right)^{m/2}\leq\infty,$$
where $\varphi_t=\varphi|_{D_t}$.
Let $E_m=\coprod_{t\in D}E_{m,t}$ be the disjoint union of all $E_{m,t}$.
Then we have a natural projection $\pi:E_m\ra D$ which maps elements in $E_{m,t}$ to $t$.
We view $H$ as a Finsler metric on $E_m$.

In general $E_m$ is not a genuine holomorphic vector bundle over $D$.
However, we can also talk about its holomorphic sections,
which are the objects we are really interested in.
By definition, a section $s:D\ra E_m$ is a \emph{holomorphic section} if it varies holomorphically with $t$,
namely, the function $s(t,z):\Omega\ra \mc$ is holomorphic with respect to the variable $t$.
Note that $s(t,z)$ is automatically holomorphic on $z$ for $t$ fixed,
by Hartogs theorem, $s(t,z)$ is holomorphic jointly on $t$ and $z$ and hence is a holomorphic function on $\Omega$.

Let $E^*_{m,t}$ be the dual space of $E_{m,t}$,
namely the space of all complex linear functions on $E_{m,t}$.
Let $E^*_m=\coprod_{t\in D}E^*_{m,t}$.
The natural projection from $E^*_m$ to $D$ is denoted by $\pi^*$.
Note that we do not define any topology on $E^*_{m,t}$ and $E_m$.
The only object we are interested in is holomorphic sections of $E^*_m$ which we are going to define.
The following definition, as well as the definition of $E_m$ given above,  is proposed in the recent work \cite{DWZZ18}.

\begin{defn}\label{def-intro:holo section of dual bundle}
A section $\xi$ of $E^*_{m}$ on $D$ is holomorphic if:
\begin{itemize}
\item[(1)] for any local holomorphic section $s$ of $E_m$, $<\xi, s>$ is a holomorphic function;
\item[(2)] for any sequence $s_j$ of holomorphic sections of $E_m$ on $D$ such that
$\int_D|s_j|_m\leq 1$, if $s_j(t,z)$ converges uniformly on compact subsets of $\Omega$ to $s(t,z)$ for some holomorphic section $s$ of $E_m$,
then $<\xi,s_j>$ converges uniformly to $<\xi,s>$ on compact subsets of $D$.
\end{itemize}
\end{defn}

In the same way we can define holomorphic section of $E^*_m$ on open subsets of $D$.
The Finsler metric $H$ on $E_m$ induces a Finsler metric $H^*$ on $E^*_m$ (see Definition \ref{def:dual finsler metric}).

\begin{thm}\label{thm-intro: psh dual section stein}
$(E_m, H)$ has positive curvature, in the sense that for any holomorphic section $\xi:D\ra E^*_m$  of $E^*_m$,
the function $\psi:=\log |\xi(t)|_m:=\log H^*(\xi(t)):D\ra [0,+\infty)$ is plurisubharmonic on $D$.
\end{thm}

For the case that $m=1$, $\Omega$ is a product, and $\varphi$ is smooth up to $\overline\Omega$
Theorem \ref{thm-intro: psh dual section stein} was proved by Berndtsson in \cite{Bob09};
and similar results with $m=1$ and varying fibers was proved by Wang \cite{Wx17} with extra technical conditions.
It is obvious from the argument that Theorem \ref{thm-intro: psh dual section stein}
can be generalized to general Stein families of complex manifolds equipped with pseudoeffective line bundles.

A direct corollary of Theorem \ref{thm-intro: psh dual section stein} is the plurisubharmonicity
of the relative $m$-Bergman kernels associated to $\Omega$ and $\varphi$,
which is proved in \cite{DWZZ18} (the case that $m=1$ is proved in \cite{Bob06}).

We now consider the case of a K\"ahler family of complex manifolds.
Let $X, Y$ be K\"ahler  manifolds of  dimension $r+n$ and $r$ respectively,
let $p:X\rightarrow Y$  be a proper holomorphic submersion.
Let $L$ be a holomorphic line bundle over $X$, and $h$ be a singular  Hermitian metric on $L$,
whose curvature current is semi-positive.

Let $m>0$ be a fixed integer.
The multiplier ideal sheaf $\mathcal I_m(h)\subset \mathcal O_X$ is defined as follows.
If $\varphi$ is a local weight of $h$ on some open set $U\subset X$,
then the germ of $\mathcal I_m(h)$ at a point $p\in U$ consists of the germs
of holomorphic functions $f$ at $p$ such that $|f|^{2/m}e^{-\varphi/m}$ is integrable at $p$.
It is known that $\mathcal I(h^{1/m})$ is a coherent analytic sheaf on $X$ \cite{Cao141}.

For $y\in Y$ let $X_y=p^{-1}(y)$ , which is a compact submanifold of $X$ of dimension $n$ if $y$ is a regular value of $p$.
Let $K_{X/Y}$ be the relative canonical bundle on $X$ and  $\mathcal E_m=p_*(mK_{X/Y}\otimes L\otimes \mathcal I_m(h))$ be the direct image sheaf on $Y$.
By Grauert's theorem, $\mathcal E_m$ is a coherent analytic sheaf on $Y$. One can choose a proper analytic subset $A\subset Y$ such that:
\begin{itemize}
\item[(1)] $p$ is submersive over $Y\backslash A$,
\item[(2)] both $\mathcal E_m$  and the quotient sheaf $p_*(mK_{X/Y}\otimes L)/\mathcal E_m$  are  locally free on $Y\backslash A$,
\item[(3)] $E_{m,y}$  is naturally identified with
$H^0(X_y,mK_{X_y}\otimes L|_{X_y}\otimes \mathcal I_m(h)|_{X_y})$, for $y\in Y\backslash A$.
\end{itemize}
where $E_m$  is the vector bundle on $Y\backslash A$
associated to $\mathcal E_m$.
For $u\in E_{m,y}$,
the $m$-norm of $u$ is defined to be
$$H_m(u):=\|u\|_m=\left(\int_{X_y}|u|^{m/2}h^{1/m}\right)^{m/2}\leq +\infty.$$
Then $H_m$ is a Finsler metric on $E_m$.

\begin{thm}\label{thm-intro:positivity Kahler case}
For any $m\geq 1$, $H_m$ is a singular Finsler metric on
the coherent analytic sheaf $\mathcal E_m$ which is positively curved.
\end{thm}

Please see Definition \ref{def:positivity of finsler} for the definition
of positively curved singular Finsler metrics on coherent analytic sheaves.

A consequence of Theorem \ref{thm-intro:positivity Kahler case} is that the $m$-Bergman kernel metric
on the twisted relative pluricanonical line bundle $mK_{X/Y}\otimes L$ has positive curvature current.
In recent years, plurisubharmonic variation of $m$-Bergman kernel metrics and
topics in the direction of Theorem \ref{thm-intro:positivity Kahler case} (in the case that $m=1$)
were extensively studied by many authors (see e.g.
\cite{Bob06}\cite{Bob09}\cite{BP08}\cite{BP10}\cite{GZh15}\cite{PT18}\cite{Cao141}\cite{HPS16}\cite{ZZ17}\cite{ZZ18}\cite{Bob18}\cite{DWZZ18}),
in different settings, in different generality, and by different methods.

Though positivity of direct image sheaves with Hermitian metrics and plurisubharmonic variation of $m$-Bergman kernels
have been extensively studied, it seems that the positivity of the Finsler metrics described in
Theorem \ref{thm-intro: psh dual section stein} and Theorem \ref{thm-intro:positivity Kahler case}
was overlooked in literatures.
Our motivation for Theorem \ref{thm-intro: psh dual section stein} and Theorem \ref{thm-intro:positivity Kahler case}
comes from Theorem \ref{thm-intro: main theorem 1} and the works of Chi and Yau as mentioned above
which imply that $E_m$ in Theorem \ref{thm-intro: psh dual section stein}
and $\mathcal E_m$ in Theorem \ref{thm-intro:positivity Kahler case}
are important linear invariants of the involved holomorphic families.

Our method to Theorem \ref{thm-intro: psh dual section stein} and Theorem \ref{thm-intro:positivity Kahler case}
follows from that in \cite{GZh15} and \cite{HPS16},
and is based on the Ohsawa-Takegoshi type extension theorems with optimal estimate (see \cite{Bl13}\cite{GZh15}\cite{Cao141}\cite{ZZ17})
and an $L^p$-variant of them \cite{BP10}.
\subsection*{Acknowledgements}
The authors are partially supported by NSFC grants.
The first author is partially supported by the University of Chinese Academy of Sciences.

\section{Linear invariants of bounded pseudoconvex domains in $\mathbb C^n$}\label{sec:linear invariants}
The aim of this section is to prove Theorem \ref{thm-intro: main theorem 1}.
The proof is divided into three subsections.
\subsection{A measure theoretic preparation}\label{subsec: Lp isometries and equimeasurability}
Let $\Omega$ be a bounded domain in $\mathbb C^n$.
Let $z=(z_1,\cdots,z_n)$ be  holomorphic coordinates of $\mathbb C^n$,
and let $d\lambda_n:={(\frac{i}{2})^n}dz_1\wedge d\overline{z}_1\wedge\cdots\wedge dz_n\wedge d\overline{z}_n$
be the canonical volume form on $\Omega$.
For $p>0$,  denote by  $A^p(\Omega)$  the space of all holomorphic functions $\phi$ with finite $L^p$-norm
\begin{align*}
\|\phi\|_p:=\left(\int_\Omega |\phi|^pd\lambda_n\right)^{1/p}.
\end{align*}
It is a standard fact that for $p\geq 1$,  $A^p(\Omega)$ are seperable Banach spaces,
and for $0<p<1$, $A^p(\Omega)$ are complete separable metric spaces with respect to the metric
$d(\varphi_1,\varphi_2):=\|\varphi_1-\varphi_2\|_p^p$.
The following result about isometries between $L^p$ spaces, is due to Rudin.

\begin{lem}[\cite{Rud76}]\label{lem:equ lp}
Let $\mu$ and $\nu$ be finite positive measures on two sets $M$ and $N$ respectively.
Assume $0<p<\infty$ and $p$ is not even.
Let $n$ be a positive integer.
If $f_i\in  L^p(M,\mu)$, $g_i\in L^p(N,\nu)$ for $1\leq i\leq n$ satisfy
\begin{align*}
\int_M|1+\alpha_1f_1+\cdots +\alpha_nf_n|^pd\mu=\int_N|1+\alpha_1g_1+\cdots +\alpha_ng_n|^pd\nu
\end{align*}
for all $(\alpha_1,\cdots,\alpha_n)\in \mathbb{C}^n$,
then $(f_1,\cdots, f_n)$ and $(g_1,\cdots, g_n)$ are equimeasurable,
i.e. for every bounded Borel measurable function (and for every real-valued nonnegative Borel function) $u:\mathbb C^n\ra \mathbb C$, we have
\begin{align*}
\int u(f_1,\cdots, f_n)d\mu=\int u(g_1,\cdots, g_n)d\nu.
\end{align*}
Furthermore, let $I:M\ra \mathbb C^n$ and $J: N\ra \mathbb C^n$ be the maps $I=(f_1,\cdots, f_n)$ and $J=(g_1,\cdots, g_n)$, respectively. Then we have
\begin{align*}
\mu(I^{-1}(E))=\nu(J^{-1}(E))
\end{align*}
for every Borel set $E$ in $\mathbb C^n$.
\end{lem}

The following lemma, observed by Markovic \cite{Mark03}, is a direct corollary of Lemma \ref{lem:equ lp}.

\begin{lem}\label{lem: equi lp higher dimension}
Let $\Omega_1$ and $\Omega_2$ be two bounded  domains in $\mathbb C^n$ and $\mathbb C^m$, repectively.
Suppose that $\phi_k$, $k=0,1,2,\cdots,N$, $N\in \mathbb N$, are elements of $A^p(\Omega_1)$ and
suppose that $\psi_k$, $k=0,1,2,\cdots, N$ are elements of $A^p(\Omega_2)$,
such that for every $N$-tuple of complex numbers $\alpha_k$, $k=1,\cdots, N$,
we have
\begin{align*}
\|\phi_0+\sum\limits_{k=1}^N\alpha_k\phi_k\|_p=\|\psi_0+\sum\limits_{k=1}^N\alpha_k\psi_k\|_p.
\end{align*}
If neither $\phi_0$ nor $\psi_0$ is constantly zero,
then for every real valued non-negative Borel function $u:\mathbb C^N\ra \mathbb C$,
we have
\begin{align*}
\int_{\Omega_1} u(\frac{\phi_1}{\phi_0},\cdots, \frac{\phi_N}{\phi_0})|\phi_0|^pd\lambda_n=\int_{\Omega_2} u(\frac{\psi_1}{\psi_0},\cdots, \frac{\phi_N}{\phi_0})|\psi_0|^pd\lambda_m.
\end{align*}
\end{lem}

\begin{proof}
Set $d\mu:=|\phi_0|^pd\lambda_n$ and $d\nu:=|\psi_0|^pd\lambda_m$.
Then the measures $\mu$ and $\nu$ are well defined on $\Omega_1$ and $\Omega_2$, respectively.
Then from the assumption, we have that for $k=1,2,\cdots,N$,
$\frac{\phi_k}{\phi_0}\in A^p(d\mu)$ and $\frac{\psi_k}{\psi_0}\in A^p(d\nu)$, respectively, and that
\begin{align*}
\int_{\Omega_1}|1+\sum\limits_{k=1}^N\alpha_k\frac{\phi_k}{\phi_0}|^pd\mu=\int_{\Omega_2}|1+\sum\limits_{k=1}^N\alpha_k\frac{\psi_k}{\psi_0}|^pd\nu.
\end{align*}
Thus Lemma \ref{lem: equi lp higher dimension} follows from Lemma \ref{lem:equ lp}.
\end{proof}

\subsection{Definition of the map and equidimension}\label{subset:define map and equidim}
Let $\Omega_1\subset\mc^n$ and $\Omega_2\subset\mc^m$ are bounded domains and let $p>0, p\neq 2,4, \cdots$.
Assume there is a linear isometry $T:A^p(\Omega_1)\ra A^p(\Omega_2)$.

We define a subset $\Omega'_1$ of $\Omega_1$ as follows.
For $z\in \Omega_1$, we say as definition that $z\in \Omega'_1$ if and only if there exits $w\in \Omega_2$ such that
\begin{equation}\label{eqn:intrinsic def map}
\phi(z)T\phi_0(w)=T\phi(w)\phi_0(z),\ \forall \phi, \phi_0\in A^p(\Omega_1).
\end{equation}
Similarly, we define $\Omega'_2\subset\Omega_2$ by saying that $w\in\Omega'_2$
if and only if $\exists z\in \Omega_1$ such that \eqref{eqn:intrinsic def map} holds.
Since $\Omega_j\ (j=1,2)$ are bounded domains, $A^p(\Omega_j)$ separates points in $\Omega_j$.
Hence for any $z\in \Omega'_1$ there is a unique $w\in \Omega_2$ such that \eqref{eqn:intrinsic def map} holds.
Therefore we can define a map
$$F:\Omega'_1\ra\Omega_2$$
by setting $F(z)=w$ if $w\in \Omega_2$ satisfies Equation  \eqref{eqn:intrinsic def map}.
It is clear that $F(\Omega'_1)\subset \Omega'_2$ and
$$F:\Omega'_1\ra\Omega'_2$$
is a bijection.

A more conceptual definition of $\Omega'_j$ and $F$ is as follows.
For $z\in \Omega_j, (j=1,2)$, $L_{j,z}:=\{\phi\in A^p(\Omega_j)|\phi(z)=0\}$ is a hyperplane in $A^p(\Omega_j)$.
For $z\in \Omega_1, w\in \Omega_2$, we say that $z\in\Omega'_1, w\in\Omega'_2$ and define $F(z)=w$ if
$T(L_{1,z})=L_{2,w}$.

From the definition, the following properties about $\Omega'_1$ and $F$ is obvious.

\begin{lem}\label{lem:coverge => in domain}
(1) For a sequence of points $z_j\in \Omega'_1$ which converges to $z_0\in \Omega_1$,
if $w_j:=F(z_j)$ converges to $w_0\in \Omega_2$, then $z_0\in\Omega'_1,\ w_0\in\Omega'_2$ and $F(z_0)=w_0$; and

(2) For $z\in\Omega'_1$ and $\phi\in A^p(\Omega)$, $\phi(z)=0$ if and only if $T\phi(F(z))=0$.
\end{lem}

To investigate more properties of $\Omega'_j\ (j=1,2)$ and $F$,
We need to give a representation of $F$.
As in \cite{Mark03}, this is down by choosing a countable dense subset of $A^p(\Omega_j)$.
Since $A^p(\Omega_1)$ is separable, we can choose a countable dense subset $\{\phi_0, \phi_1, \cdots\}$ of it.
Then $\{\psi_0:=T\phi_0, \psi_1:=T\phi_1, \cdots\}$ is a countable dense subset of $A^p(\Omega_2)$.
For $N>0$, we define maps $I_N, J_N$ as follows:
$$I_N:\Omega_1\ra\mc^N,\ z\mapsto (\frac{\phi_1(z)}{\phi_0(z)},\cdots, \frac{\phi_N(z)}{\phi_0(z)}),$$
$$J_N:\Omega_2\ra\mc^N,\ z\mapsto (\frac{\psi_1(z)}{\psi_0(z)},\cdots, \frac{\psi_N(z)}{\psi_0(z)}),$$
which are measurable maps on $\Omega_1$ and $\Omega_2$ respectively.
We can also define $I_\infty$ and $J_\infty$ as
$$I_\infty:\Omega_1\ra\mc^\infty,\ z\mapsto (\frac{\phi_1(z)}{\phi_0(z)}, \frac{\phi_2(z)}{\phi_0(z)}, \cdots),$$
$$J_\infty:\Omega_2\ra\mc^\infty,\ z\mapsto (\frac{\psi_1(z)}{\psi_0(z)}, \frac{\psi_2(z)}{\psi_0(z)}, \cdots).$$
Then it is clear that $\Omega'_1\backslash\phi^{-1}_0(0)=\cap_N I^{-1}_N J_N(\Omega_2)=I_\infty^{-1}J_\infty(\Omega_2)$
and $\Omega'_2\backslash\psi^{-1}_0(0)=\cap_N J_N^{-1}I_N(\Omega_1)=J_\infty^{-1}I_\infty(\Omega_1)$.
Since $A^p(\Omega_j)$ separates points of $\Omega_j, \ (j=1,2)$,
both $I_\infty$ and $J_\infty$ are injective on their domains of definition.
For $z\in I_\infty^{-1}J_\infty(\Omega_2)$, it is easy to see that $F(z)=J^{-1}_\infty(I_\infty(z))$.

\begin{lem}\label{lem: D' full measure}
The Lebesgue measure of $\Omega_j\backslash\Omega'_j$ are zero, for $j=1, 2$.
\end{lem}
\begin{proof}
It suffices to prove that Lebesgue measure of $K:=\Omega_1\backslash\Omega'_1$ is zero.
Set $d\mu:=|\phi_0|^pd\lambda_n$ and $d\nu:=|\psi_0|^pd\lambda_m$.
It suffices to prove that $\mu(K)=0$.

By Lemma \ref{lem: equi lp higher dimension}, we have
$\mu(I^{-1}_N(I_N(K)))=\nu(J^{-1}_N(I_N(K)))$.
We have seen that $\Omega'_1\supseteq\cap_N I^{-1}_N J_N(\Omega_2)=I_\infty^{-1}J_\infty(\Omega_2)$,
so $\cap_N J_N^{-1}I_N(\Omega_1\backslash\Omega'_1)=J_\infty^{-1}I_\infty(\Omega_1\backslash\Omega'_1)=\emptyset$.
Note also that $J^{-1}_N(I_N(K))$ is decreasing with respect to $N$,
and $\nu(\Omega_2)<\infty$, we have $\lim_N \nu(J^{-1}_N(I_N(K)))=0$.
Since $\mu(K)\leq \mu(I^{-1}_N(I_N(K)))$ for all $N$, we see $\mu(K)=0$.
\end{proof}

\begin{prop}\label{prop:equi dimension}
The dimensions of $\Omega_1$ and $\Omega_2$ are equal, namely $n=m$.
\end{prop}
\begin{proof}
We first prove that $n\geq m$.
Choose $\phi_0, \phi_1, \cdots$ such that $\psi_0=1, \psi_1=w_1, \cdots, \psi_m=w_m$,
where $(w_1,\cdots, w_m)$ are coordinates on $\mc^m$.
It is clear that $J_\infty^{-1}(I_\infty(\Omega_1))\subset J_m^{-1}(I_m(\Omega_1))$.
By Lemma \ref{lem: D' full measure} (and its proof), $J_\infty^{-1}(I_\infty(\Omega_1))$ has full measure in $\Omega_2$,
hence $J_m^{-1}(I_m(\Omega_1))$ has full measure in $\Omega_2$.
By the choice of $\phi_0, \phi_1, \cdots$, $J_m$ is just the identity map,
and hence $I_m(\Omega_1)$ has positive measure in $\mc^m$.
Note that $I_m$ is a holomorphic map on $\Omega_1\backslash \phi^{-1}_0(0)$, we have $n\geq m$.
The same argument shows that $m\geq n$.
So we get $m=n$.
\end{proof}

\begin{lem}\label{lem: D' open}
$\Omega'_j$ are open subsets of $\Omega_j$, for $j=1, 2$.
\end{lem}
\begin{proof}
Let $z_0\in\Omega'_1$ and $w_0\in \Omega'_2$ such that $F(z_0)=w_0$.
Choose $\phi_0, \phi_1, \cdots$ such that $\psi_0=1, \psi_1=w_1, \cdots, \psi_m=w_m$.
Then $J_m:\Omega_2\ra\mc^n$ is the identity map.

By Lemma \ref{lem:coverge => in domain} (2), we have $\phi_0(z_0)\neq 0$.
So $I_m$ is holomorphic on some neighborhood $U$ of $z_0$ in $\Omega_1$.
By contracting $U$ if necessary, we may assume $I_m(U)\subset\Omega_2$.

It is obvious that $F(z)=I_m(z)$ for $z\in U\cap\Omega'_1$,
which implies that
$$\phi(z)T\phi'(I_m(z))=T\phi(I_m(z))\phi'(z), \ \forall \phi, \phi'\in A^p(\Omega_1), \forall z\in U\cap\Omega'_1.$$
By Lemma \ref{lem: D' full measure}, $U\cap\Omega'_1$ is dense in $U$.
By continuity, the above equality holds for all $z\in U$.
By definition of $\Omega'_1$ and $F$, we have $U\subset \Omega'_1$ and $F=I_m$ on $U$.
This implies $\Omega'_1$ is open.
By the same argument, one can prove that $\Omega'_2$ is open.
\end{proof}

\begin{lem}\label{lem:T1=Jacobi}
For $z\in\Omega'_1$ and $\phi\in A^p(\Omega_1)$, we have
$$|T\phi(F(z))||J_F(z)|^{2/p}=|\phi(z)|, \ \phi\in A^p(\Omega_1), z\in\Omega'_1,$$
where $J_F$ is the holomorphic Jacobian of $F$.
\end{lem}
\begin{proof}
Let $K$ be any measurable subset in $\Omega_1$.
Considering the function $u(t_1, \cdots t_N)=|t_1|$ on $\mc^N$,
from Lemma \ref{lem: equi lp higher dimension}, we have
$$\int_{I_N^{-1}(I_N(K))}|\phi_1|^pd\lambda_n=\int_{J_N^{-1}(I_N(K))}|\psi_1|^pd\lambda_n, \ \forall N.$$
Note that $I_N^{-1}(I_N(K))$ decrease to $K$ and $J_N^{-1}(I_N(K))$ decrease to $F(K)$ as $N\ra\infty$,
we get
$$\int_{K}|\phi_1|^pd\lambda_n=\int_{F(K)}|\psi_1|^p d\lambda_n,.$$
Similarly the above equality holds for all $\phi_i, i\geq 1$.
Considering that $\{\phi_0,\phi_1,\cdots\}$ is dense in $A^p(\Omega_1)$,
we get
$$\int_{K}|\phi|^pd\lambda_n=\int_{F(K)}|T\phi|^p d\lambda_n,\ \forall \phi\in A^p(\Omega_1).$$
On the other hand, we have
$$\int_{F(K)}|T\phi|^p d\lambda_n=\int_{K}|T\phi(F(z))|^p|J_F(z)|^2.$$
So we get
$$\int_{K}|\phi(z)|^pd\lambda_n=\int_{K}|T\phi(F(z))|^p|J_F(z)|^2$$
for all Borel subset $K$ in $\Omega_1$.
By continuity, we have
$$|\phi(z)|^p=|T\phi(F(z))|^p|J_F(z)|^2,\ z\in\Omega'_1.$$
\end{proof}

\subsection{Extending $F$ to a biholomorphic map between $\Omega_1$ and $\Omega_2$.}\label{subsect: biholomorphism construction}
We complete the proof of Theorem \ref{thm: main theorem 1}.
We follow the notations in \ref{subset:define map and equidim}.
For convenience, we restate Theorem \ref{thm: main theorem 1} here:

\begin{thm}\label{thm: main theorem 1}
Let $\Omega_1\subset\mc^n$ and $\Omega_2\subset\mc^m$ be bounded hyperconvex domains.
Suppose that there is a $p>0, p\neq 2,4,6 \cdots$ , such that
\begin{itemize}
\item[(1)] there is a  linear isometry $T:A^p(\Omega_1)\ra A^p(\Omega_2)$,  and
\item[(2)] the $p$-Bergman kernels of  $\Omega_1$ and $\Omega_2$ are exhaustive,
\end{itemize}
then $m=n$ and there exists a  biholomorphic map $F:\Omega_1\ra \Omega_2$ such that
\begin{align*}
|T\phi\circ F||J_F|^{2/p}=|\phi|, \ \forall \phi\in A^p(\Omega_1),
\end{align*}
where $J_F$ is the holomorphic Jacobian of $F$.
If $n=1$, the assumption of hyperconvexity can be dropped.
\end{thm}

\begin{proof}
We show that $F:\Omega'_1\ra \Omega'_2$ can extend to a biholomorphic map from $\Omega_1$ to $\Omega_2$.
From Proposition \ref{prop:equi dimension}, we have $n=m$.
The rest of the proof is divided into several lemmas.

\begin{lem}\label{lem: definition of A}
Let $S=\Omega_1\backslash\Omega'_1$, then for any $\zeta\in S$, $\lim\limits_{\Omega'_1\ni z\rightarrow \zeta}J_{F}(z)=0$.
\end{lem}

\begin{proof}
We argue by contradiction.
Suppose to the contrary, there exists $\{z_j\}\in \Omega'_1$,
such that $z_j\rightarrow \zeta$ and $|J_F(z_j)|\geq \varepsilon>0$ for all $j$.
Since $\Omega_2$ is bounded, we may also assume that $F(z_j)\rightarrow w\in \mathbb C^n$.
Then either $w\in \Omega_2$ or $w\in \partial\Omega_2$.

If $w\in \Omega_2$,
by Lemma \ref{lem:coverge => in domain}, $\zeta\in \Omega'_1$.
This is a contradiction.

We now have $w\in \partial\Omega_2$.
Since the $p$-Bergman kernel of $\Omega_2$ is exhaustive,
we have that  $\lim\limits_{j\rightarrow +\infty}K_{\Omega_2,p}(F(z_j))\ra +\infty$.
Thus there are  $\psi_j\in A^p(\Omega_2)$, such that $\|\psi_j\|_p=1$ and $|\psi_j(F(z_j))|\rightarrow +\infty$.
Let $\psi_j=T(\varphi_j)$, then $\|\varphi_j\|_p=1$.
By Lemma \ref{lem:T1=Jacobi}, we have  $ |(\psi_j\circ F)||{ J_F}|^{2/p}=|\varphi_j|$ on $\Omega'_1$.
Thus $|\varphi_j(z_j)|=|\psi_j(F(z_j))||J_F(z_j)|^{2/p}\ra +\infty$.
This means  that  $K_{\Omega_1,p}(z_j)\ra +\infty$, which is a contradiction,
since $z_j\ra \zeta\in \Omega_1$ and  $K_{\Omega_1,p}$ is locally bounded on $\Omega_1$.
\end{proof}

We define a function $\rho:\Omega_1\rightarrow \mr\cup\{-\infty\}$ as
$\rho(z)=\log|J_F(z)|$ for $z\in \Omega'_1$ and $\rho(z)=-\infty$ as $z\in\Omega_1\backslash \Omega'_1$.

\begin{lem}\label{lem: psh of rho}
The function $\rho$ is a plurisubharmonic function on $\Omega_1$ and $\rho^{-1}(-\infty)=\Omega_1\backslash \Omega'_1$.
\end{lem}

\begin{proof}
Since $f$ is biholomorphic on $\Omega'_1$, $\rho(z)\neq -\infty$ for $z\in \Omega'_1$.
By definition of $\rho$, we have $\rho^{-1}(-\infty)=\Omega_1\backslash \Omega'_1$.
From Lemma \ref{lem: definition of A}, $\rho$ is upper semicontinuous.
Since $\rho$ is plurisubharmonic on $\Omega'_1$ and, by Lemma \ref{lem: D' open}, $\Omega_1\backslash \Omega'_1$ is closed,
it is clear that $\rho$ satisfies the mean value inequality when restricted on
any complex line intersection with $\Omega_1$.
Hence $\rho$ is plurisubharmonic.
\end{proof}

Recall that a subset $A$ in a complex manifold $X$ is called a pluripolar set if
for any $z\in X$ there is a neighborhood $U$ of $z$ and a plurisubharmonic function
$h$ on $U$ such that $A\cap U\subset h^{-1}(-\infty)$.
Lemma \ref{lem: psh of rho} implies that $S:=\Omega_1\backslash \Omega'_1$
is a closed pluripolar set of $\Omega_1$.

We need the following well-known fact about continuation of holomorphic functions.
\begin{lem}[see e.g. \cite{Demailly}]\label{lem: extension from closed pluripolar set}
Let $X$ be complex analytic manifold and let $A\subset X$ be a closed pluripolar set.
Then\\
(1) every holomorphic function on $X\setminus A$ that is locally bounded near any point in $A$
extends to a holomorphic function on $X$; and \\
(2) every plurisubharmonic function on $X\setminus A$ that is locally bounded above near any point in $A$
extends to a plurisubharmonic function on $X$.
\end{lem}

We have constructed a biholomorphic map $F: \Omega_1\setminus S\ra \Omega_2$.
From Lemma \ref{lem: extension from closed pluripolar set},
$F$ extends across $S$ to a holomorphic map, also denoted by $F$,  from $\Omega_1$ to $\mc^n$.
We want to show that $S=\emptyset$ and $\Omega'_1=\Omega_1$.

We argue by contradiction.
If $S\neq \emptyset$, by Lemma \ref{lem:coverge => in domain} and note that $\Omega_1\backslash S$ is dense in $\Omega_1$,
$F(S)\subset \overline\Omega_2\backslash\Omega_2$.

Since $\Omega_2$ is hyperconvex, there is a plurisubharmonic function $\varphi$ on $\Omega_2$
such that $\varphi<0$ and for any $c<0$ the set $\{w\in\Omega_2| \varphi(w)\leq c\}$ is compact.

Let $\tilde\varphi=\varphi\circ F$, then $\tilde\varphi$ is a plurisubharmonic function on $\Omega_1\backslash S$
which have 0 as an upper bound.
By Lemma \ref{lem: extension from closed pluripolar set}, $\tilde\varphi$ extends to a plurisubharmonic function
on $\Omega_1$ which attains its maximum on $S$.
By the maximum principle of plurisubharmonic functions, $\tilde\varphi$ is constant.
This is a contradiction.
So $S=\emptyset$ and $\Omega'_1=\Omega_1$.

By exchanging the roles of $\Omega_1$ and $\Omega_2$,
we can prove that $\Omega'_2=\Omega_2$.
Hence $F:\Omega_1\ra\Omega_2$ is a biholomorphic map.

For the case that $n=1$, it is clear that $S=\{z\in\Omega_1| F'(z)=0\}$
and hence is a discrete subset of $\Omega_1$.
Note that $F:\Omega'_1\ra\Omega_2$ is injective, $F'(z)\neq 0$ for all $z\in \Omega_1$
and hence $S=\emptyset$ and $\Omega'_1=\Omega_1$.
Similarly, we have $\Omega'_2=\Omega_2$.
Hence $F$ is a biholomorphic map from $\Omega_1$ to $\Omega_2$.

The equality in Theorem \ref{thm: main theorem 1} (2) follows from Lemma \ref{lem:T1=Jacobi}.

We now prove the uniqueness of $F$.
If there is another biholomorphic map $G:\Omega_1\ra\Omega_2$
satisfying the conditions in Theorem \ref{thm: main theorem 1}.
Then $G$ satisfies
$$\phi(z)T\phi_0(G(z))=T\phi(G(z))\phi_0(z),\ \forall z\in\Omega_1, \forall \phi\in A^p(\Omega_2).$$
By the definition of $F$ (see \S \ref{subset:define map and equidim}), we have $G=F$.
\end{proof}

\section{Exhaustion of $p$-Bergman kernels}
Theorem \ref{thm: main theorem 1} implies
the importance of the exhaustion property of $p$-Bergman kernels for bounded domains.
In this Section, we discuss exhaustion of $p$-Bergman kernels in various settings.

For $0<p<2$, Ning-Zhang-Zhou proved the following complete result.

\begin{thm}[\cite{NZZ16}]\label{thm: exhaustion}
For any bounded domain $\Omega$ in $\mathbb C^n$ and any $p\in (0,2)$,
$\Omega$ is pseudoconvex if and only if the $p$-Bergman kernel $B_{\Omega,p}(z)$ is exhaustive.
\end{thm}

\begin{proof}
For the sake of completeness, we include the proof here.

Since $B_{\Omega,p}$ is a plurisubharmonic function on $\Omega$,
$\Omega$ is pseudoconvex provided $B_{\Omega,p}$ is exhaustive.

We now prove that $B_{\Omega,p}$ is exhaustive if $\Omega$ is pesudoconvex.
For $\zeta\in\partial D$, we need to show that $\lim_{z\ra \zeta}B_{\Omega, p}(z)=+\infty$.
After a unitary transform,
we may assume that $\zeta=(a,0,\cdots, 0)$.
Let $L$ be the line given by $L=\{z_2=\cdots=z_n=0\}$.
It is clear that $f(z_1):=\frac{1}{z_1-a}\in A^p(\Omega\cap L)$.
By a $L^p$ variant of the Ohsawa-Takegoshi extension theorem obtained in \cite{BP10},
there is $F\in A^p(\Omega)$ such that $F|_{\Omega\cap L}=f$,
which implies that $\lim_{z\ra \zeta}B_{\Omega, p}(z)=+\infty$.
\end{proof}

Recall that a bounded domain $\Omega\subset \mc^n$ is called homogeneous regular
if there is a constant $c\in (0,1)$ such that for any $z\in \Omega$,
there is a holomorphic injective map $f:\Omega\ra \mathbb B^n$ with
$f(z)=0$ and $\mathbb B^n(c)\subset f(\Omega)$,
where $\mathbb B^n(c)=\{z\in\mc^n; \|z\|<c\}$.

\begin{thm}\label{thm: exhaustion simply connected USD}
If $\Omega\subset\mc^n$ is a simply connected and homogeneous regular domain,
then the $p$-Bergman kernel of $\Omega$ is exhaustive.
 \end{thm}

 To prove Theorem \ref{thm: exhaustion simply connected USD}, we need two lemmas.
 The first one is about the transformation formula for $p$-Bergman kernels.

 \begin{lem}[\cite{NZZ16}]\label{lem: transition of Bergman kernels}
 Let $\Omega_1, \Omega_1$ be simply connected domains in $\mathbb C^n$.
 Then for any biholomorphism  $f:\Omega_1\ra \Omega_2$,
 we have $B_{\Omega_1,p}(z)=B_{\Omega_2,p}(f(z))|J_F|^2$,
 where $J_F$ is the holomorphic Jacobian of $f$.
 Furthermore, if $p=\frac{2}{m}$, where $m\in \mathbb N$,
 there is no need to assume that the domains are simply connected.
 \end{lem}

 \begin{lem}[Generalized Rouch\' e's theorem, c.f. \cite{Llo79}] \label{lem: generalized Rouche}
 Let $\Omega$ be a bounded domain of $\mathbb C^n$.
 Let $f$ and $g$ be two continuous mappings of $\overline{\Omega}$ (the closure of $\Omega$)
 into $\mathbb C^n$ which are holomorphic in $\Omega$.
 Suppose that
 \begin{align*}\|g(z)\|<\|f(z) \|~~~~~~~~~(z\in \partial\Omega)
 \end{align*}
 Then $f$ and $f+g$ has the same number of zeros in $\Omega$, counting multiplicities,
 where $\|\cdot\|$ is the standard norm in $\mathbb C^n$.
 \end{lem}

 \begin{proof}[Proof of Theorem  \ref{thm: exhaustion simply connected USD}]
 By assumption, there is a positive constant $c>0$, such that $\forall z\in \Omega$,
 there is a holomorphic injection $f_z:\Omega\ra \mathbb B^n$,
 such that $f_z(z)=0$, $\mathbb B^n(c)\subset f_z(\Omega)$.

 Fix an arbitrary point $a\in \partial \Omega$.
 For any sequence of points  $z_j\in \Omega$ that converges to $a$,
 we need to prove that
 \begin{align*}
 \lim\limits_{j\ra \infty}B_{\Omega,p}(z_j)=+\infty.
 \end{align*}

We argue by contradiction.
Suppose to the contrary,  we may assume that $B_{\Omega,p}(z_j)\leq M$ for all $j$, for some positive constant $M$.
For simplicity we denote by $f_j$ the map $f_{z_j}$.
Since $\Omega$ is simply connected, from Lemma \ref{lem: transition of Bergman kernels}, we have
 \begin{align*}
 B_{\Omega,p}(z_j)=B_{f_j(\Omega,p)}(0)|J_{f_j}(z_j)|^2.
 \end{align*}
From the decreasing property of $p$-Bergman kernels, we have
 \begin{align*}
 B_{f_j(\Omega),p}(0)\geq B_{\mathbb B^n,p}(0).
 \end{align*}
 Thus $M\geq B_{\Omega,p}(z_j)\geq B_{\mathbb B^n,p}(0)|J_{f_j}(z_j)|^2$,
 which implies that $|J_{f_j}(z_j)|\leq M'$ for all $j$, for some constant $M'>0$.
 By Montel Theorem, we may assume that $f_j\ra f$,  and $g_j:=f_j^{-1}|_{\mathbb B^n(c)}\ra g:\mathbb B^n(c)\ra \mathbb C^n$.
 Then $|J_g(0)|=\lim\limits_j |J_{g_j}(0)|\geq M'^{-1}$.
  In particular, $J_g(0)\neq 0$, and hence $g$ is a locally biholomorphic mapping at $0$ with $g(0)=\lim\limits_j g_j(0)=\lim\limits_j z_j=a\in \partial \Omega$.
  Now choose a sufficiently small neighborhood $U\Subset \mathbb B^n(c)$ of $0$, such that $g|_{\partial U}\neq a$.
 Then map $g-a$ has a zero inside $U$.
 From Lemma \ref{lem: generalized Rouche}, $g_j-a$ has a zero inside $U$ for all $j$ sufficiently large,
 which is a  contradiction since $g_j\neq a$.
\end{proof}

For bounded strongly pseudoconvex domains with $\mathcal C^2$ boundary, which are homogeneous regular \cite{Deng-Guan-Zhang16},
the same statement in Theorem \ref{thm: exhaustion simply connected USD} holds,
with the simply connectedness of $\Omega$ dropped.
The proof is based on results on exposing boundary points of strongly pseudoconvex domains in \cite{Diedrich-Fornaess-Wold13}\cite{DFW18}.

\begin{thm}\label{thm: exhaustion strongly pscd}
For any bounded strongly pseudoconvex domain $\Omega$ in $\mathbb C^n$
with $\mathcal C^2$ smooth boundary, the $p$-Bergman kernel is exhaustive for all $p>0$.
\end{thm}
\begin{proof}
Let $a$ be an arbitrary boundary point of $\Omega$.
Let $R>0$ be a sufficiently large positive number, such that $\overline\Omega\subset\mathbb B^n(R)=\{z\in\mc^n; \|z\|<R\}$.
Then for any (small) neighborhood $U$ of $a$ in $\mc^n$ and any $\epsilon>0$,
there is a holomorphic injective map $f$ defined in some neighborhood of $\overline\Omega$
such that $f(a)=(R,0,\cdots, 0)$, $f(\Omega)\subset\mathbb B^n(R)$, and $\|f(z)-z\|<\epsilon$ for $z\in\Omega\backslash U$ (see \cite{DFW18}\cite{Diedrich-Fornaess-Wold13}).

Let $f(z)=(f_1(z), \cdots, f_n(z))$.
Taking $U$ small enough, we can have $f_1\in A^p(\Omega)$ with $\|f_1\|_p\leq M$,
where $M>0$ is a constant that can be chosen to be independent of $R$.
Letting $R\ra +\infty$, we get $B_{\Omega, p}(z)\ra +\infty$ as $z\ra a$.
\end{proof}

\section{Singular Finsler metrics on coherent analytic sheaves}

In this section,
we recall the notions introduced in \cite{DWZZ18} of singular Finsler metrics on holomorphic vector bundles
and give a definition of positively curved singular Finsler metrics on coherent analytic sheaves.

\begin{defn}\label{def:finsler on v.b.}
Let $E\rightarrow X$ be a holomorphic vector bundle over a complex manifold $X$. A (singular) Finsler metric $h$ on $E$
is a function $h:E\ra [0,+\infty]$, such that $|v|^2_h:=h(cv)=|c|^2h(v)$ for any $v\in E$ and $c\in\mc$.
\end{defn}
In the above definition, we do not assume any regularity property of $h$.
Only when considering Griffiths positivity certain regularity is required,
as shown in the following Definition \ref{def:finsler on v.b.}.

\begin{defn}\label{def:dual finsler metric}
For a singular Finsler metric $h$ on $E$,
its dual Finsler metric $h^*$ on the dual bundle $E^*$ of $E$ is defined as follows.
For $f\in E^*_x$, the fiber of $E^*$ at $x\in X$,
$|f|_{h^*}$ is defined to be
$$|f|_{h^*}:=\sup\{|f(v)|; v\in E_x, |v|_h\leq 1\}\leq +\infty.$$
\end{defn}

\begin{defn}\label{def:positivity of finsler}
Let $E\rightarrow X$ be a holomorphic vector bundle over a complex manifold $X$.
A singular Finsler metric $h$ on $E$ is called negatively curved (in the sense of Griffiths)
if for any local holomorphic section $s$ of $E$ the function $\log|s|^2_h$ is plurisubharmonic,
and is called positively curved (in the sense of Griffiths) if its dual metric $h^*$ on $E^*$ is negatively curved.
\end{defn}

Let $\mathcal F$ be a coherent analytic sheaf on $X$,
it is well known that $\mathcal F$ is locally free on some Zariski open set $U$ of $X$.
On $U$, we will identify $\mathcal F$ with the vector bundle associated to it.

\begin{defn}\label{def:finsler on sheaf}
Let $\mathcal{F}$ be a coherent analytic sheaf on a complex manifold $X$.
Let $Z\subset X$ be an analytic subset of $X$ such that $\mathcal{F}|_{X\setminus Z}$ is locally free.	
A negatively curved singular Finsler metric $h$ on $\mathcal{F}$ is a singular Finsler metric on the holomorphic vector bundle $\mathcal{F}|_{X\setminus Z}$,
such that for any local holomorphic section $s$ of $\mathcal{F}$ on any open set $U\subset X$,
the function $\log|s|_{h}$ is plurisubharmonic on  $U\setminus Z$, and can be extended to a plurisubharmonic function on $U$;
and a singular Finsler metric $h$ on a coherent analytic sheaf $\mathcal F$
is said to be positively curved if the dual metric $h^*$ on $\mathcal{F}^*$ is negatively curved.
\end{defn}

\begin{rem}\label{rem: extension of finsler metric}
Suppose that $\log|g|_{h^*}$  is p.s.h. on $U\setminus Z$.
It is well-known that if  codim$_{\mathbb{C}}(Z)\geq 2$ or $\log|g|_{H^*}$ is locally bounded above near $Z$,
then $\log|g|_{h^*}$ extends across $Z$ to $U$ uniquely as a p.s.h function.
Definition \ref{def:finsler on sheaf} matches
Definition \ref{def:positivity of finsler} if $\mathcal F$ is a vector bundle.
\end{rem}

\section{Variations of linear invariants-for families of pseudoconvex domains}
In this section, we will prove Theorem \ref{thm-intro: psh dual section stein}.

Let $\Omega\subset \mathbb{C}^{r+n}=\mathbb{C}_t^r\times \mathbb{C}_z^n$ be a pseudoconvex domain.
Let $p:\Omega\rightarrow \mathbb{C}^r$ be the natural  projection.
We denote $p(\Omega)$ by $D$ and denote $p^{-1}(t)$ by $\Omega_t$ for $t\in D$.
Let $\varphi$ be a plurisubharmonic function on $\Omega$ and let $m\geq 1$ be a fixed integer.
For an open subset $U$ of $D$, we denote by $\mathcal F(U)$ the space of holomorphic functions $F$ on $p^{-1}(U)$ such that
$\int_{p^{-1}(K)}|F|^{2/m}e^{-\varphi}\leq\infty$ for all compact subset $K$ of $D$.
For $t\in D$, let
$$E_{m,t}=\{F|_{\Omega_t}:F\in\mathcal F(U),\ U\subset D\ \text{open\ and}\ U\ni t\}.$$
Then $E_{m,t}$ is a vector space with the following pseudonorm:
$$H(f):=|f|_m=\left(\int_{D_t}|f|^{2/m}e^{-\varphi_t}\right)^{m/2}\leq\infty,$$
where $\varphi_t=\varphi|_{D_t}$.
Let $E_m=\coprod_{t\in D}E_{m,t}$ be the disjoint union of all $E_{m,t}$.
Then we have a natural projection $\pi:E_m\ra D$ which maps elements in $E_{m,t}$ to $t$.
We view $H$ as a Finsler metric on $E_m$.

By definition, a section $s:D\ra E_m$ is a \emph{holomorphic section} if it varies holomorphically with $t$,
namely, the function $s(t,z):\Omega\ra \mc$ is holomorphic with respect to the variable $t$.
Note that $s(t,z)$ is automatically holomorphic on $z$ for $t$ fixed,
by Hartogs theorem, $s(t,z)$ is holomorphic jointly on $t$ and $z$ and hence is a holomorphic function on $\Omega$.

Let $E^*_{m,t}$ be the dual space of $E_{m,t}$,
namely the space of all complex linear functions on $E_{m,t}$.
Let $E^*_m=\coprod_{t\in D}E^*_{m,t}$.
The natural projection from $E^*_m$ to $D$ is denoted by $\pi^*$.
Note that we do not define any topology on $E^*_{m,t}$ and $E_m$.
The only object we are interested in is holomorphic sections of $E^*_m$ which we are going to define.
The following definition, as well as the definition of $E_m$ given above,  is proposed in the recent work \cite{DWZZ18}.

\begin{defn}\label{def:holo section of dual bundle}
A section $\xi$ of $E^*_{m}$ on $D$ is holomorphic if:
\begin{itemize}
\item[(1)] for any local holomorphic section $s$ of $E_m$, $<\xi, s>$ is a holomorphic function;
\item[(2)] for any sequence $s_j$ of holomorphic sections of $E_m$ on $D$ such that
$\int_D|s_j|_m\leq 1$, if $s_j(t,z)$ converges uniformly on compact subsets of $\Omega$ to $s(t,z)$ for some holomorphic section $s$ of $E_m$,
then $<\xi,s_j>$ converges uniformly to $<\xi,s>$ on compact subsets of $D$.
\end{itemize}
\end{defn}

In the same way we can define holomorphic section of $E^*_m$ on open subsets of $D$.
The Finsler metric $H$ on $E_m$ induces a Finsler metric $H^*$ on $E^*_m$ (defined in the same way as in Definition \ref{def:dual finsler metric}).

The metric $H$ and $H^*$ have the following semicontinuity properties
which are proved in \cite{DWZZ18}.

\begin{prop}\label{prop:semi-continuous: noncompact}
Assume $s$ is a holomorphic section of $E_m$,
then the function $|s|_m(t):=H(s(t)):D\ra (0,+\infty]$ is lower semicontinuous;
and let $\xi:D\ra E^*_m$ be a holomorphic section of $E^*_m$.
Then the function $|\xi|_m(t):=H^*(\xi(t)):D\ra [0,+\infty)$ is upper semicontinuous.
\end{prop}

The following result is the main ingredient in our proof or Theorem \ref{thm-intro: psh dual section stein}.

\begin{lem}[\cite{DWZZ18}]\label{lem:Lm-extension on psc}
Let $\Omega\subset \mathbb{C}^{n+1}$ be a pseudoconvex domain, $\Delta$ be the unit disk in $\mathbb{C}$, and
$p: \Omega\rightarrow p(\Omega)\subset \Delta$ be  a holomorphic  projection.
For $y\in \Delta$, we denote $\Omega_y:=p^{-1}(y)$ by $\Omega_y$.
Let $\varphi$  be a p.s.h function on $\Omega$. Let $m\geq 1$ be an integer and $y_0\in \Delta$ such that $\varphi$ is not identically $-\infty$ on any branch of $\Omega_y$.
Then for any holomorphic function $u$ on $\Omega_{y_0}$ such that
\begin{align*}
\int_{\Omega_{y_0}}|u|^{2/m}e^{-\varphi}<+\infty,
\end{align*}
there exists a holomorphic function $U$ on $\Omega$ such that $U|_{\Omega_{y_0}}=u$ and
\begin{align*}
\int_{\Omega}|U|^{2/m}e^{-\varphi}\leq \pi \int_{\Omega_{y_0}}|u|^{2/m}e^{-\varphi},
\end{align*}
\end{lem}
The proof of Lemma \ref{lem:Lm-extension on psc} is a combination of the iteration method of Berndtsson-P\u aun \cite{BP10}
and the optimal Ohsawa-Takegoshi extension theorems (\cite{Bl13} \cite{GZh15}).

We can now give the proof of Theorem \ref{thm-intro: psh dual section stein}
which we restate it here.

\begin{thm}\label{thm: psh dual section stein}
$(E_m, H)$ has positive curvature, in the sense that for any holomorphic section $\xi:D\ra E^*_m$  of $E^*_m$,
the function $\psi:=\log |\xi(t)|_m:=\log H^*(\xi(t)):D\ra [0,+\infty)$ is plurisubharmonic on $D$.
\end{thm}
\begin{proof}
Form Proposition \ref{prop:semi-continuous: noncompact},
we already know that $\psi$ is upper semi continuous on $D$.
Thus it suffices to prove that for every holomorphic mapping $\gamma:\Delta\rightarrow D$ from the unit disc $\Delta$ to $D$,
the function $\psi$ satisfies the mean-value inequality
\begin{align*}
(\psi\circ \gamma)(0)\leq \frac{1}{\pi}\int_{\Delta}(\psi\circ\gamma).
\end{align*}
Without loss of generality, we may assume that $D=\Delta$ and $\gamma$ is the identity map.
The above inequality holds trivially if $\psi(0)=-\infty$.
We now assume that $\psi(0)\neq -\infty$.
Then we can choose an element $f\in E_{m,0}$ such that $|f|_{H,0}=1$ and
\begin{align*}
\psi(0)=\log |\xi|_{H,0}=\log |\xi(f)|.
\end{align*}
By Lemma \ref{lem:Lm-extension on psc},  there is a holomorphic section $s\in H^0(\Delta, E_{m})$, such that
\begin{align*}
s(0)=f, ~~~\frac{1}{\pi}\int_{\Delta}|s|_{H}^{2/m}d\mu\leq 1
\end{align*}
By definition of the metric $H^*$ on $E_m^*$, we have the pointwise inequality
\begin{align*}
|\xi|_{H^*}\geq\frac{|\xi(s)|}{|s|_{H}}
\end{align*}
and therefore, $\psi\geq \log|\xi(s)|-\log |s|_{H}$. Multiplying $2/m$ both sides, we get that
\begin{align*}
\frac{2}{m} \psi\geq\frac{2}{m}\log|\xi(s)|-\log |s|_{H}^{2/m}
\end{align*}
Integrating, we get
\begin{align*}
\frac{1}{\pi}\int_{\Delta} \frac{2}{m}\psi d\mu\geq \frac{1}{\pi}\int_{\Delta}\frac{2}{m}\log|\xi(s)|d\mu-\frac{1}{\pi}\int_{\Delta}\log|s|_{H}^{2/m}d\mu
\end{align*}

Note that  $\log|\xi(s)|$ is p.s.h., thus satisfies the mean-value inequality, and so the first term on the R.H.S.  is at least $\frac{2}{m}\log|\xi(f)|=\frac{2}{m}\psi(0)$. Since the function $|s|^{2/m}_{H}$ is integrable, from Jensen's inequality we get that
\begin{align*}
-\frac{1}{\pi}\int_{\Delta}\log|s|_{H}^{2/m}d\mu \geq -\log(\frac{1}{\pi}\int_{\Delta}|s|^{2/m}_{H}d\mu)\geq -\log 1=0.
\end{align*}
In conclude, we obtain that
\begin{align*}
\frac{1}{\pi}\int_{\Delta} \frac{2}{m}\psi d\mu\geq \frac{2}{m}\psi(0),
\end{align*}
which means that $\frac{2}{m}\psi$ is subharmonic, i.e. $\psi$ is subharmonic.
\end{proof}


A consequence of Theorem \ref{thm: psh dual section stein} is the following
\begin{cor}\label{cor: psh for holomorphic measure}
For $t\in D$, let $t\rightarrow \mu_t$ be a family of complex measures on $\Omega$.
Assume that $\mu_t$ has compact support along fibers in the sense that for any $t_0\in D$,
there is an open subset $U\subset D$ with $t_0\in U$,
such that $supp\mu_t\subset K\ (t\in U)$ for some compact subset $K\subset\Omega$. Then
\begin{align*}u_t\ra \int_\Omega u_td\mu_t=:\mu_t(u_t)
\end{align*}
defines a section of the dual bundle $E_m^*$.
Suppose that the section $\mu_t$ is holomorphic in the sense that
\begin{align*}
t\ra \mu_t((t,\cdot))
\end{align*}
is holomorphic for any local holomorphic section $s$ of $E_m$.
Then $\log\|\mu_t\|_{H^*}$ is plurisubharmonic on $D$.
\end{cor}

When the measures $\mu_t$ are all Dirac delta functions, we get the following
\begin{cor}[\cite{DWZZ18}]\label{cor: psh for variation of m-Bergman kernels}
For any $m\geq 1$,
the function $\log K_{m,t}(z)$ is a plurisubharmonic function on $\Omega$,
where $K_{m,t}(z)$ is the relative twisted $2/m$-Bergman kernel on $\Omega$ with weight $\varphi$.
\end{cor}
\begin{rem}
The case of $m=1, n=1, \varphi=0$ is proved by Mataini and Yamaguchi \cite{MY04},
and that of $m=1$ and general $n$ and $\varphi$ is proved by Berndtsson \cite{Bob06}.
The general form is  proved in \cite{DWZZ18} by applying a new characterization of p.s.h.
functions and  the Ohsawa-Takegoshi extension theorem.
\end{rem}

\section{Variations of linear invariants-for families of K\"ahler manifolds}\label{sec:E_m positive Kaler case}
The aim of this section is to prove Theorem \ref{thm-intro:positivity Kahler case}.

We first recall the definition of canonical pseudonorms on the space of twisted pluricanonical sections,
and their semicontinuity property associated to families of compact K\"ahler manifolds.

Let $X, Y$ be K\"ahler  manifolds of  dimension $r+n$ and $r$ respectively,
let $p:X\rightarrow Y$  be a proper holomorphic submersion.
Let $L$ be a holomorphic line bundle over $X$, and $h$ be a singular  Hermitian metric on $L$,
whose curvature current is semi-positive.

Let $m>0$ be a fixed integer.
The multiplier ideal sheaf $\mathcal I_m(h)\subset \mathcal O_X$ is defined as follows.
If $\varphi$ is a local weight of $h$ on some open set $U\subset X$,
then the germ of $\mathcal I_m(h)$ at a point $p\in U$ consists of the germs
of holomorphic functions $f$ at $p$ such that $|f|^{2/m}e^{-\varphi/m}$ is integrable at $p$.
It is known that $\mathcal I(h^{1/m})$ is a coherent analytic sheaf on $X$ \cite{Cao141}.

For $y\in Y$ let $X_y=p^{-1}(y)$ , which is a compact submanifold of $X$ of dimension $n$ if $y$ is a regular value of $p$.
Let $K_{X/Y}$ be the relative canonical bundle on $X$ and  $\mathcal E_k=p_*(mK_{X/Y}\otimes L\otimes \mathcal I_m(h))$ be the direct image sheaf on $Y$.
By Grauert's theorem, $\mathcal E_m$ is a coherent analytic sheaf on $Y$. One can choose a proper analytic subset $A\subset Y$ such that:
\begin{itemize}
\item[(1)] $p$ is submersive over $Y\backslash A$,
\item[(2)] both $\mathcal E_m$  and the quotient sheaf $p_*(mK_{X/Y}\otimes L)/\mathcal E_m$  are  locally free on $Y\backslash A$,
\item[(3)] $E_{m,y}$  is naturally identified with
$H^0(X_y,mK_{X_y}\otimes L|_{X_y}\otimes \mathcal I_m(h)|_{X_y})$, for $y\in Y\backslash A$.
\end{itemize}
where $E_m$  is the vector bundle on $Y\backslash A$
associated to $\mathcal E_m$.
For $u\in E_{m,y}$,
the $m$-norm of $u$ is defined to be
$$H_m(u):=\|u\|_m=\left(\int_{X_y}|u|^{m/2}h^{1/m}\right)^{m/2}\leq +\infty.$$
Then $H_m$ is a Finsler metric on $E_m$.

About the semicontinuity of the metrics $H_m$ and $H^*_m$, we have
\begin{prop}[\cite{HPS16}]\label{prop:Hodge metric semi-continuous:compact}
Let $s$ be a holomorphic section of $E_m$.
The function $|s|_m(y):=\|s(y)\|_m:Y\rightarrow (0,+\infty]$ is lower semi-continuous.
For every $\xi\in H^0(Y, E^*_m)$, the function $|\xi|_{H_m^*}(y):=H_m^*(\xi(y)):Y\ra [0,+\infty)$ is  upper semi-continuous.
\end{prop}

For the proof of Theorem \ref{thm-intro:positivity Kahler case}, we need the following

\begin{thm}[\cite{DWZZ18}]\label{lem:Lm-extension projective} Let $B\subset \mc^r$ be the unit ball and let $X$ be a K\"ahler manifold of dimension $n+r$.
Let $p:X\ra B$ be a holomorphic proper submersion.
For $t\in B$, denote by $X_t$ the fiber $p^{-1}(t)$.
Let $L$ be a holomorphic line bundle on $X$ with a (singular) Hermitian metric $h$
whose curvature current is positive. Let $u\in H^0(X_0, mK_{X_0}\otimes \mathcal I_m(h)|_{X_0})$ with $\|u\|_m<\infty$.
Assume there exists an open subset $U$ containing the origin and $s_0\in H^0(U,\mathcal E_m|_U)$ such that $s_0|_{X_0}=u$,
then there exists $s\in H^0(X,mK_X\otimes L\otimes \mathcal I_m(h))=H^0(B, mK_B\otimes \mathcal E_m)$
such that $s|_{X_0}=u\wedge dt^{\otimes m}$ and
$$\int_X |s|^{2/m}h^{1/m}\leq \mu(B)\int_{X_0}|u|^{m/2}h^{1/m},$$
where $t=(t_1,\cdots, t_r)$ is the standard coordinate on $B$ and $dt=dt_1\wedge\cdots\wedge dt_r$,
and $\mu(B)$ is the volume of $B$ with respect to the Lebesgue measure on $B$.
\end{thm}
The proof of this theorem is a combination of the iteration method of Berndtsson-P\u aun \cite{BP10}
and the optimal Ohsawa-Takegoshi extension theorems on weakly pseudoconvex K\"{a}hler manifolds (see \cite{ZZ}).

We can now give the proof of Theorem \ref{thm-intro:positivity Kahler case}.
For convenience, we restate it here.

\begin{thm}\label{thm:positivity Kahler case}
For any $m\geq 1$, $H_m$ is a singular Finsler metric on
the coherent analytic sheaf $\mathcal E_m$ which is positively curved.
\end{thm}

\begin{proof}
The argument is similar to that in the proof of Theorem \ref{thm-intro: psh dual section stein}.

From Definition \ref{def:finsler on sheaf}, we need to prove that,
for every $g\in H^0(Y,\mathcal{E}_m^*) $, the funtion $\psi=\log |g|_{H_m^*}$ is p.s.h. on $Y\setminus A$,
and can be extended to be a p.s.h. function on $Y$.

From Proposition \ref{prop:Hodge metric semi-continuous:compact},
the function $\psi$ is upper semi-continuous on $Y\setminus A$.
From Remark \ref{rem: extension of finsler metric},
it suffices to prove that $\psi$ is p.s.h.  on $Y\setminus A$ and locally bounded above near $A$.

We first prove the plurisubharmonicity of $\psi$ on $Y\setminus A$.
It suffices to prove that for every holomorphic mapping $\gamma:\Delta\rightarrow Y\setminus A$,
such that $\psi\circ\gamma$ is not identically $-\infty$,
$\psi$ satisfies the mean-value inequality
\begin{align*}
(\psi\circ \gamma)(0)\leq \frac{1}{\pi}\int_{\Delta}(\psi\circ\gamma)d\mu.
\end{align*}

We may assume that $Y\backslash A=\Delta$ and assume that $\psi(0)\neq 0$.
Then there exits an element $f\in E_{m,0}$ with $|f|_{H_m}=1$ and
\begin{align*}
\psi(0)=\log |g|_{H_m,0}=\log |g(f)|.
\end{align*}
By Cartan theorem B and Lemma \ref{lem:Lm-extension projective},
there is a holomorphic section $s\in H^0(\Delta, E_{m})$, such that
\begin{align*}
s(0)=f, ~~~\frac{1}{\pi}\int_{\Delta}|s|_{H_m}^{2/m}d\mu\leq 1
\end{align*}
By definition of the metric $H^*_m$ on the dual bundle, we have the pointwise inequality
\begin{align*}
|g|_{H^*_m}\geq\frac{|g(s)|}{|s|_{H_m}}
\end{align*}
and therefore, $\psi\geq \log|g(s)|-\log |s|_{H_m}$. Multiplying $2/m$ both sides, we get that
\begin{align*}
\frac{2}{m} \psi\geq\frac{2}{m}\log|g(s)|-\log |s|_{H_m}^{2/m}
\end{align*}
Integrating over $\Delta$, we get
\begin{align*}
\frac{1}{\pi}\int_{\Delta} \frac{2}{m}\psi d\mu\geq \frac{1}{\pi}\int_{\Delta}\frac{2}{m}\log|g(s)|d\mu-\frac{1}{\pi}\int_{\Delta}\log|s|_{H_m}^{2/m}d\mu.
\end{align*}

Since $g(s)$ is a holomorphic function on $\Delta$,
$\log|g(s)|$ is p.s.h. and hence satisfies the mean-value inequality.
It thus suffice to prove that the second term in the above inequality is $\leq$ 0.
This is true from Jensen's inequality which implies
$$\frac{1}{\pi}\int_{\Delta}\log|s|_{H_m}^{2/m}d\mu\leq \log(\frac{1}{\pi}\int_{\Delta}|s|^{2/m}_{H_m}d\mu)\leq \log 1=0.$$
In conclude, we obtain that
$\frac{2}{m}\psi$ is subharmonic, so $\psi$ is subharmonic.

We now give the proof of the boundeddness of $\psi$.
Without loss of generality, we assume that $Y=\mathbb B^r$ is the unit ball.
For any $y\in Y\backslash A$ such that $|g(y)|\neq 0$ (otherwise there is nothing to prove),
there is $a\in E_{m,y}$ such that $|a|=1$ and $\langle g(y),a\rangle=|g(y)|$.
By Cartan theorem B and  Lemma \ref{lem:Lm-extension projective},
there exists a holomorphic section $s$ of $\mathcal E_m$ on $Y$
such that $s(y)=a$ and
$$\int_{Y}|s|^{2/m}_{H_m}\leq C,$$
where $C$ is a constant independent of $y$ and $a$.
Let
$$S=\{f\in H^0(Y,\mathcal E_m); \int_{Y}|f|^{2/m}_{H_m}\leq C\}.$$
Since the metric $h$ on $L$ is lower semicontinuous,
by the mean value inequality and Montel theorem, $S$ is a normal family,
namely, any sequence in $S$ has a subsequence that converges uniformly on compact subsets of $Y$.
Note that if $s_j$ is a subsequence of $S$ that converges uniformly on compact subsets of $Y$,
then the sequence of holomorphic functions $\langle u,s_j\rangle$ converges on compact subsets of $Y$,
and hence is uniformly bounded on compact sets of $Y$.
So $\{\langle u,s\rangle;s\in S\}$ is uniformly bounded on compact sets of $Y$.

We no complete the proof of the theorem \ref{thm:positivity Kahler case}.
\end{proof}

We can consider the pull back $\tilde {E}_m:=p^* E_m$,
which is a coherent analytic sheaf on $X\backslash p^{-1}(A)$.
The Finsler metric $H_m$ on $E_m$ naturally induces a metric $\tilde H_m$ on $\tilde{E}_m$
which is clearly positively curved in the sense of Definition \ref{def:finsler on sheaf}.
The evaluation map from $\tilde{ E}_m|_{X\backslash p^{-1}(A)}$ to $(mK_{X/Y}\otimes L)|_{X\backslash p^{-1}(A)}$
induces a hermitian metric (the quotient metric) on $(mK_{X/Y}\otimes L)|_{X\backslash p^{-1}(A)}$,
which is called the $m$-Bergman kernel metric on $(mK_{X/Y}\otimes L)|_{X\backslash p^{-1}(A)}$.
From Theorem \ref{thm:positivity Kahler case}, we get

\begin{cor}
The relative $m$-Bergman kernel metric defined above on $(mK_{X/Y}\otimes L)|_{X\backslash p^{-1}(A)}$ has nonnegative curvature current,
and can be extended to a singular Hermitian metric on $mK_{X/Y}\otimes L$ whose curvature current is nonnegative.
\end{cor}

\subsection{Application to Teichm\"{u}ller metric}
In this section, we apply Theorem \ref{thm:positivity Kahler case} to show that
the Teichm\"uller metric on Teichm\"{u}ller spaces are negatively curved.
For basicn knowledge about Teichm\"{u}ller spaces please see \cite{IT92} and \cite{St84}.

Let $C$ be a compact Riemann surface of genus $g\geq 2$.
Given a quasiconformal mapping $f:C\rightarrow S$,
the pair $(S,f)$ defined to be a marked Riemannian surface.
Two marked Riemannian surfaces $(S_1,f_1)$ and $(S_2,f_2)$ are called equivalent
if $f_2\circ f_1^{-1}: S_1\ra S_2$ is homotopic to a conformal mapping $c:S_1\ra S_2$.
Denote by $[S,f]$ the equivalence class of $(S,f)$.
The Teichm\"{u}ller space $T=T(C)$ is the set of equivalence classes of all marked Riemann surfaces $[S,f]$.

It is well-known that $T$ carries a natural complex structure
and there is a holomorphic family $p:X\ra T$ of compact Riemann surfaces
such that for each $t\in T$ the fiber $X_t:=p^{-1}(t)$ belongs to the
equivalence class that is represented by $t$.

The cotangent space $T^*_tT$ of $T$ at a $t$ can be identified with the space $H^0(t,2K_{X_t})$
of holomorphic quadratic differentials on $X_t$.
There is a natural Finsler metric on $T$ which is called the Teichm\"uller metric on $T$.
The norm on $T^*_tT$ dual to the norm on $T_tT$ associated to the Teichm\"uller metric is given by
\begin{align*}
\|s\|_1:=\int_S(s\wedge \overline s)^{1/2}.
\end{align*}

Now we consider the coheret sheaves $\mathcal E_m:=p_*(mK_{X/T})$ equipped with
the canonical Finsler metric $H_m$ introduced in Section \ref{sec:E_m positive Kaler case}.
As explained above, $\mathcal E_2$ can be identified with the sheaf associated to the cotangent bundle of $T$
and $H_2$ is dual to the Techm\"uller metric on the tangent bundle of $T$.

As a direct corollary of Theorem \ref{thm:positivity Kahler case},
we have the following
\begin{cor}\label{cor: m positivity teichmuller}
$H_m$ is a positively curved Finsler metric on $\mathcal E_m$ .
\end{cor}

When $m=2$, we have

\begin{thm}\label{thm: negative tangent}
The Teichm\"{u}ller metric on $T$ is a negatively curved Finsler metric.
\end{thm}


\begin{thebibliography}{10}
\bibitem{Ave-Hed-Per15}B. Avelin, L. Hed, H. Persson,
A note on the hyperconvexity of pseudoconvex domains beyond Lipschitz regularity. (English summary)
Potential Anal. 43 (2015), no. 3, 531¨C545.

\bibitem{Bob18}
B.~Berndtsson.
\newblock Complex Brunn-Minkowski theory and positivity of vector bundles.
\newblock {\em arXiv:1807.05844}.

\bibitem{Bob06}
B.~Berndtsson.
\newblock Subharmonicity properties of the {B}ergman kernel and some other
  functions associated to pseudoconvex domains.
\newblock {\em Ann. Inst. Fourier (Grenoble)}, 56(6):1633--1662, 2006.

 \bibitem{Bob09}
    B.~Berndtsson.
    \newblock Curvature of vector bundles associated to holomorphic fibrations.

   \bibitem{BP08}
    B.~Berndtsson and M.~P\u{a}un.
    \newblock Bergman kernels and the pseudoeffectivity of relative canonical
    bundles.
    \newblock {\em Duke Math. J.}, 145(2):341--378, 2008.    \newblock {\em Ann. of Math. (2)}, 169(2):531--560, 2009.

\bibitem{BP10}
B.~Berndtsson and M.~P\u{a}un.
\newblock Bergman kernels and subadjunction.
\newblock {\em arXiv:1002.4145}.

\bibitem{Bl13}
Z.~B\l ocki.
\newblock Suita conjecture and the {O}hsawa-{T}akegoshi extension theorem.
\newblock {\em Invent. Math.}, 193(1):149--158, 2013.

\bibitem{Cao141}
J.~Cao.
\newblock Ohsawa-{T}akegoshi extension theorem for compact {K}\"ahler manifolds
  and applications.
\newblock In {\em Complex and symplectic geometry}, volume~21 of {\em Springer
  INdAM Ser.}, pages 19--38. Springer, Cham, 2017.

\bibitem{Chen-Zhang}B. Chen, J. Zhang, \emph{On graphs of holomorphic motions}, preprint.
\bibitem{Chi2016} C.-Y. Chi, Pseudonorms and theorems of Torelli type, J. Differential Geom. 104 (2016), no. 2, 239-273.

\bibitem{Chi-Yau2008} C.-Y. Chi and S.-T. Yau, A geometric approach to problems in birational geometry. Proc. Natl. Acad. Sci. USA 105 (2008), no. 48, 18696-18701.

\bibitem{Demailly87}J.P. Demailly, Mesures de Monge-Amp\'ere et mesures pluriharmoniques, Math. Z. (1987) vol. 194, 519-564.
\bibitem{Demailly}J.P. Demailly, Complex analytic and differential geometry, e-book, available at:http://www-fourier.ujf-grenoble.fr/~demailly/documents.html.
\bibitem{DFW18} F. Deng, J. E. Forn{\ae}ss, E. F. Wold, Exposing boundary points of strongly pseudoconvex subvarieties in complex
spaces, Proc. Amer. Math. Soc. 146(2018), 2473-2487.
\bibitem{Deng-Guan-Zhang12}F. Deng, Q. Guan, L. Zhang, \emph{On some properties of squeezing functions on bounded domains}, Pacific J. Math. V. 57, No.2 (2012), 319-342.
\bibitem{Deng-Guan-Zhang16}F. Deng, Q. Guan, L. Zhang, \emph{Properties of squeezing functions and global transformations of bounded domains}, Transactions of AMS, 368 (2016), 2679-2696.

\bibitem{DWZZ18}
F.~Deng, Z.~Wang, L.~Zhang, and X.~Zhou.
\newblock New characterization of plurisubharmonic functions and positivity of
  direct image sheaves.
\newblock {\em Preprint}, 2018.
\bibitem{Diedrich-Fornaess-Wold13}K. Diederich, J.E. Forn{\ae}ss, E.F. Wold, \emph{Exposing points on the boundary of a strictly pseudoconvex or a locally convexifiable domain of finite 1-type }, Journal of Geometric Analysis,  24 (2014), 2124-2134.
\bibitem{Fornaess-Wold16-preprint}J.E. Forn{\ae}ss, E.F. Wold, \emph{A non-strictly pseudoconvex domain for which the squeezing function tends to one towards the boundary}, e-preprint, arXiv:1611.04464.
\bibitem{Frankel91}S. Frankel, Applications of affine geometry to geometric function theory in several complex variables, I: Convergent rescalings and
intrinsic quasi-isometric structure, Several Complex Variables and Complex Geometry, Part 2 (Santa Cruz, CA, 1989), 183-208, Proc.
YSympos. Pure Math., 52, Part 2, Amer. Math. Soc., Providence, RI, 1991.
\bibitem{GZh15}
Q.~Guan and X.~Zhou.
\newblock A solution of an {$L^2$} extension problem with an optimal estimate
  and applications.
\newblock {\em Ann. of Math. (2)}, 181(3):1139--1208, 2015.

\bibitem{Lakic97}N. Lakic,
An isometry theorem for quadratic differentials on Riemann surfaces of finite genus.
Trans. Amer. Math. Soc. 349 (1997), no. 7, 2951¨C2967.

\bibitem{HPS16}
C.~Hacon, M.~Popa, and C.~Schnell.
\newblock Algebraic fiber spaces over abelian varieties: {A}round a recent
  theorem by {C}ao and {P}\u aun.
\newblock In {\em Local and global methods in algebraic geometry}, volume 712
  of {\em Contemp. Math.}, pages 143--195. Amer. Math. Soc., Providence, RI,
  2018.

\bibitem{IT92}
Y.~Imayoshi and M.~Taniguchi.
\newblock {\em An introduction to {T}eichm\"uller spaces}.
\newblock Springer-Verlag, Tokyo, 1992.
\newblock Translated and revised from the Japanese by the authors.
\bibitem{Kim-Zhang16} K.-T. Kim, L. Zhang, \emph{On the uniform squeezing property and the squeezing function}, Pacific J. Math. 282 (2016), 341-358.
\bibitem{LS14}
L.~Lempert and R.~Sz\H oke.
\newblock Direct images, fields of {H}ilbert spaces, and geometric
  quantization.
\newblock {\em Comm. Math. Phys.}, 327(1):49--99, 2014.

\bibitem{Liu-Sun-Yau04}K. Liu, X. Sun, S.T. Yau, \emph{Canonical metrics on the moduli space of Riemann Surfaces I}, J. Differential Geom. Vol. 68 (2004), 571-637.

\bibitem{Llo79}
N.~G. Lloyd.
\newblock Remarks on generalising {R}ouch\'{e}'s theorem.
\newblock {\em J. London Math. Soc. (2)}, 20(2):259--272, 1979.

\bibitem{MY04}
F.~Maitani and H.~Yamaguchi.
\newblock Variation of {B}ergman metrics on {R}iemann surfaces.
\newblock {\em Math. Ann.}, 330(3):477--489, 2004.

\bibitem{Mark03}
V.~Markovic.
\newblock Biholomorphic maps between {T}eichm\"{u}ller spaces.
\newblock {\em Duke Math. J.}, 120(2):405--431, 2003.

\bibitem{Nikolov-Andreev16}N. Nikolov, L. Andreev, \emph{Boundary behavior of the squeezing functions of complex domains}, e-preprint, arXiv:1609.02051.
\bibitem{NZZ16}
J.~Ning, H.~Zhang, and X.~Zhou.
\newblock On {$p$}-{B}ergman kernel for bounded domains in {$\Bbb C^n$}.
\newblock {\em Comm. Anal. Geom.}, 24(4):887--900, 2016.

\bibitem{PT18}
M.~P\u aun and S.~Takayama.
\newblock Positivity of twisted relative pluricanonical bundles and their
  direct images.
\newblock {\em J. Algebraic Geom.}, 27(2):211--272, 2018.

\bibitem{Rud76}
W.~Rudin.
\newblock {$L^{p}$}-isometries and equimeasurability.
\newblock {\em Indiana Univ. Math. J.}, 25(3):215--228, 1976.

\bibitem{St84}
K.~Strebel.
\newblock {\em Quadratic differentials}, volume~5 of {\em Ergebnisse der
  Mathematik und ihrer Grenzgebiete (3) [Results in Mathematics and Related
  Areas (3)]}.
\newblock Springer-Verlag, Berlin, 1984.

\bibitem{Wx17}
X.~Wang.
\newblock A curvature formula associated to a family of pseudoconvex domains.
\newblock {\em Ann. Inst. Fourier (Grenoble)}, 67(1):269--313, 2017.
  \bibitem{Yau15}
    S.-T. Yau.
    \newblock On the pseudonorm project of birational classification of algebraic
    varieties.
    \newblock In {\em Geometry and analysis on manifolds}, volume 308 of {\em
        Progr. Math.}, pages 327--339. Birkh\"auser/Springer, Cham, 2015.
\bibitem{Yeung09}S. K. Yeung, \emph{Geometry of domains with the uniform squeezing property}, Adv.  Math. 221 (2009) 547-569.


\bibitem{ZZ17}
X.~Zhou and L.~Zhu.
\newblock An optimal $L^2$ extension theorem on weakly pseudoconvex k\"{a}hler
manifolds.
\newblock {\em J. of Differential Geom}, 110(1):135--186, 2018.

\bibitem{ZZ18}
X.~Zhou and L.~Zhu.
\newblock Siu's lemma, optimal $L^2$ extension, and applications to
pluricanonical sheaves.
\newblock {\em to appear in Math. Ann.}
\end{thebibliography}
\end{document}